\documentclass[11pt,a4paper,sunlimits]{article}
\usepackage{amsmath,amsfonts}
\usepackage{amsthm}
\usepackage{amssymb}
\usepackage{amscd}
\usepackage{xypic}
\usepackage{verbatim}
\usepackage[plainpages=false,colorlinks,hyperindex,pdfpagemode=None,bookmarksopen,linkcolor=red,citecolor=blue,urlcolor=blue]{hyperref}
\usepackage{pdflscape}
\usepackage{stmaryrd}
\usepackage{mathabx}
\usepackage{multirow}
\usepackage{appendix}
 \swapnumbers
\theoremstyle{plain}

\renewcommand{\marginpar}[1] {  }
\renewcommand{\comment}[1] {  }

\usepackage[all]{xy}
\catcode`\é=\active\def é{\'e}
\catcode`\è=\active\def è{\`e}
\catcode`\à=\active\def à{\`a}
\catcode`\ù=\active\def ù{\`u}

\catcode`\ê=\active\def ê{\^{e}}
\catcode`\î=\active\def î{\^{i}}
\catcode`\ô=\active\def ô{\^{o}}
\catcode`\û=\active\def û{\^{u}}
\catcode`\ç=\active\def ç{\c{c}}

\newtheorem{theo}{Theorem}[section]
\newtheorem{lem}[theo]{Lemma}

\theoremstyle{remark}

\numberwithin{equation}{section}

\def\bb{\backslash}
\def\ga{\gamma}

\def\la{\lambda}

\def\si{\sigma}

\def\bb{\backslash}

\def\cB{{\mathcal B}}

\def\cL{{\mathcal L}}
\def\cM{{\mathcal M}}

\def\cW{{\mathcal W}}

\DeclareFontFamily{OT1}{rsfs}{}
\DeclareFontShape{OT1}{rsfs}{n}{it}{<-> rsfs10}{}
\DeclareMathAlphabet{\mathscr}{OT1}{rsfs}{n}{it}

\def\E{{\rm E}}%{\mathbb E}
\def\F{{\rm F}}%{\mathbb F}

\def\C{\mathbb C}

\def\sG{\underline{ G}}
\def\sH{\underline{ H}}

\def\sS{\underline{ S}}

\def\me{\medskip}

\def\ste{\par\smallskip\noindent}
\def\ste{\par\smallskip\noindent}

\def\dem{ {\em Proof~: \ste }}
 \def\beq{\begin{equation}}
\def\eeq{\end{equation}}

\newenvironment{res}
               {\begin{equation}\begin{minipage}{0.85\textwidth}}
               {\end{minipage}\end{equation}}
\def\ber{\begin{res}}
\def\eer{\end{res}}

\def\qed{{\null\hfill\ \raise3pt\hbox{\framebox[0.1in]{}}\break\null}}

\begin{document}

\author{P. Delorme\thanks{The first author was supported by a grant of Agence Nationale de la Recherche with reference
ANR-13-BS01-0012 FERPLAY.}, P. Harinck}
\title{Spherical character of a supercuspidal representation as weighted orbital integral}
\date{}
\maketitle

\begin{abstract} Let $\E/\F$ be an unramified  quadratic extension of local non archimedean  fields of characteristic $0$. Let $\sH$ be an algebraic reductive group, defined and split over $\F$. We assume that the split connected component of the center of $\sH$ is trivial. Let $(\tau,V)$ be a $\sH(\F)$-distinguished supercuspidal representation of $\sH(\E)$. Using the recent results of C. Zhang \cite{Z}, and  the geometric side of a local relative trace formula obtained by P.~Delorme, P. Harinck and S.~Souaifi \cite{DHSo}, we describe spherical characters associated to $\sH(\F)$-invariant linear forms on $V$ in terms of weighted orbital integrals of matrix coefficients of $\tau$.\end{abstract}
\noindent{\it Mathematics Subject Classification 2000:} 11F72, 22E50. \medskip

\noindent{\it Keywords and phrases:} $p$-adic reductive groups, symmetric spaces, truncated kernel,  spherical character, weighted orbital integral.
\section{Introduction}

Let $\E/\F$ be an unramified  quadratic extension of local non archimedean  fields of characteristic $0$. Let $\sH$ be an algebraic reductive group, defined and split over $\F$. We denote by $\sG:=\textrm{Res}_{\E/\F}\sH_{/\E}$ the restriction of scalars of $\sH_{/\E}$. Then $G:=\sG(\F)$ is isomorphic to $\sH(\E)$. We set $H:=\sH(\F)$. We denote by $\si$ the involution of $\sG$ induced by the nontrivial element of the Galois group of $\E/\F$.

 An unitary irreducible admissible representation $(\pi,V )$ of $G$ is $H$-distinguished if the space $V^{*H}=\textrm{Hom}_H(\pi,\C)$ of $H$-invariant linear forms on $V$ is nonzero. In that case, a  distribution $m_{\xi,\xi'}$, called spherical character, can be associated to  two $H$-invariant linear forms $\xi,\xi'$ on $V$ (cf. (\ref{defcoef})).  By (\cite{Ha} Theorem 1), spherical characters are  locally integrable functions on $G$, which are $H$ biinvariant and smooth  on the set $G^{\si-reg}$ of elements $g$, called $\si$-regular points, such that $g$ is semisimple and $g^{-1}\si(g)$ is regular in $G$ in the usual sense.  \me

%From now, we assume that $A_H=\{1\}$.  
% In \cite{Z},  gives a description of $V^{*H}$ when   $(\pi, V)$ is a $H$-distinguished supercuspidal representation of $G$.

We assume that the split component of the center of $H$ is trivial.  Let $(\tau, V)$ be a $H$-distinguished supercuspidal representation of $G$. 

The aim of this note is to give the value of a spherical character $m_{\xi,\xi'}(g)$, when $g\in G$   is a regular point for the symmetric space  $H\bb G$ and $\xi,\xi'\in V^{*H}$, in terms of weighted orbital integrals of a matrix coefficient of $\tau$ (cf. Theorem \ref{coefWOI} ). This result is analogous to that of J. Arthur in the group case (\cite{ArWOI}). Notice that this result of J.~Arthur can be deduced from  his local trace formula (\cite{ArLT}) which was obtained later.\me

We use the recent work  of   C.~ Zhang \cite{Z}, which describes the space of $H$-invariant linear forms of supercuspidal representations, and the geometric side of a local relative trace formula obtained by P.~Delorme, P. Harinck and S.~Souaifi \cite{DHSo}.
\section{Spherical characters}
We denote by 
$C_c^\infty(G)$ the space of compactly smooth functions on $G$.\\
We fix a $H$-distinguished supercuspidal  representation $(\tau,V)$ of $G$.  We denote by $d(\tau)$ its formal degree.\me

Let $(\cdot,\cdot)$ be a $G$-invariant hermitian inner product on $V$. Since $\tau$ is unitary, it induces an isomorphism $\iota:v\mapsto (\cdot, v)$  from  the conjugate complex vector space $\overline{V}$ of $V$ and the smooth dual $\check{V}$ of $V$, which intertwines the complex conjugate of 
$\tau$ and its contragredient $\check{\tau}$. If $\xi$ is a linear form on $V$, we define the linear form $\overline{\xi}$ on $\overline{V}$ by $\bar{\xi}(u):=\overline{\xi(u)}.$\\
For $\xi_1$ and $\xi_2$  two  $H$-invariant linear forms on $V$, we associate the spherical character $m_{\xi_1,\xi_2}$ defined to be  the distribution on $G$  given by
\beq\label{defcoef}m_{\xi_1,\xi_2}(f):=\sum_{u\in \cB} \xi_1\big( \tau(f) u\big)\overline{\xi_2(u)},\eeq
where $\cB$ is an orthonormal basis of $V$. Since $\tau(f)$ is of finite rank, this sum is finite. Moreover, this sum does not depend on the choice of $\cB$. Indeed,  let $(\tau^*,V^*)$ be the dual representation of $\tau$. For $f\in C_c^\infty(G)$, we set $\check{f}(g):=f(g^{-1})$.   By (\cite{R} Théorème III.3.4 and I.1.2), the linear form $\tau^*(\check{f})\xi$ belongs to $\check{V}$. Hence we can write   $\iota^{-1}(\tau^*(\check{f})\xi)=\sum_{v\in \cB}\big(\tau^*(\check{f})\xi\big)(v)\cdot v$  where $(\la,v)\mapsto \la\cdot v$ is the action of $\C$ on $\overline{V}$. Therefore we deduce easily that one has 
\beq\label{coefgen}m_{\xi_1,\xi_2}(f)=\overline{\xi}_2\big( \iota^{-1}(\tau^*(\check{f})\xi_1)).\eeq\\
Since $\tau$ is a supercuspidal representation, we can define the $H\times H$-invariant pairing $\cL$ on $V\times \overline{V}$ by
$$\cL(u,v):=\int_{H} (\tau(h)u,v) dh.$$
By (\cite{Z} Theorem 1.5), 
\ber\label{surj}the map $v\mapsto \xi_v:u\mapsto \cL(u,v)$ is a surjective linear map from $\overline{V}$ onto $V^{*H}$.\eer\\
For $v,w\in V$, we denote by $c_{v,w}$ the corresponding matrix coefficient  defined by  $c_{v,w}(g):=(\tau(g)v,w)$ for $g\in G$. 
%%%%%%%%
\begin{lem}\label{calculcoef}Let $\xi_1,\xi_2\in V^{*H}$ and $v,w\in V$. Then we have
$$m_{\xi_1,\xi_2}(\check{c}_{v,w})= d(\tau)^{-1}\xi_1(v) \overline{\xi_2(w)}.$$
\end{lem}
\dem By (\ref{surj}), there exist $v_1$ and $v_2$ in $V$ such that $\xi_j=\xi_{v_j}$ for $j=1,2$.  By definition of the spherical character, for $f\in C_c^\infty(G)$ and $\cB$ an orthonormal basis of $V$, one has

$$m_{\xi_1,\xi_2}(f)=\sum_{u\in\cB}\int_{H}(\tau(h)\tau(f)u, v_1)dh \int_{H} \overline{(\tau(h)u, v_2)} dh$$
$$=\sum_{u\in\cB}\int_{H \times H }(u,\tau(\check{f})\tau(h_1) v_1)(\tau(h_2)v_2,u) dh_1 dh_2$$
$$=\int_{H\times H}(\tau(h_2)v_2,\tau(\check{f})\tau(h_1) v_1) dh_1 dh_2$$
Hence we obtain
\beq\label{express1}m_{\xi_1,\xi_2}(f)=\int_{H\times H}\int_Gf(g) (\tau(h_1gh_2)v_2,v_1) dg dh_1 dh_2.\eeq\\
Let  $f(g):=\check{c}_{v,w}(g)= \overline{(\tau(g)w,v)}$. By the orthogonality relation of Schur, for $h_1,h_2\in H$, one has
$$\int_G (\tau(g)\tau(h_2)v_2,\tau(h_1)v_1)\overline{(\tau(g)w,v)} dg= d(\tau)^{-1} (\tau(h_2)v_2, w) (v,\tau(h_1)v_1).$$
Thus, we deduce that 
$$m_{\xi_1,\xi_2}( f)=  d(\tau)^{-1} \xi_w(v_2) \xi_{v_1}(v)=d(\tau)^{-1}\xi_1(v) \overline{\xi_2(w)}.\quad \quad\qed$$\me

\section{Main result}
We first recall some notations of \cite{DHSo} to introduce weighted orbital integrals. 

We refer the reader to 
(\cite{RR} \textsection 3) and (\cite{DHSo} \textsection 1.2 and  1.3) for the notations below and more details on $\si$-regular points. Let $D_G$ be the usual Weyl discriminant function of $G$. By (\cite{RR} Lemma 3.2 and Lemma 3.3), an element $g\in G$ is $\si$-regular if and only if $D_G(g^{-1}\si(g))\neq 0$. The set  $G^{\si-reg}$  of $\si$-regular points of $G$ is decribed    as follows. Let $\sS$ be a maximal torus of $\sH$. We denote by $\sS_\si$ the connected component of the set of points $\ga \in \textrm{Res}_{\E/\F}\sS_{/\E}$ such that $\si(\ga)=\ga^{-1}$. We set $S_\si:=\sS_\si(\F)$. By Galois cohomology, there exists a finite set $\kappa_S\subset G$ such that $\sH\sS_\si\cap G=\cup_{x\in\kappa_S} HxS_\si$.\\
By (\cite{RR} Theorem 3.4) and (\cite{DHSo} (1.30)), if   $g\in G^{\si-reg}$, there exist a unique maximal torus $\sS$ of $\sH$ defined over $\F$ and 2 unique points $x\in\kappa_S$ and $\ga\in  S_\si$ such that $g=x\ga$. We denote by $M$ the centralizer of the  split connected component of $S:=\sS(\F)$. Then $M$ is Levi subgroup, that is the Levi component of a parabolic subgroup of $H$. We define the weight function $w_M$ on $H\times H$ by
$$w_M(y_1,y_2):=\tilde{v}_M(1,y_1,1,y_2),$$
where $\tilde{v}_M$ is the weight function defined in (\cite{DHSo} Lemma 2.10)
and $1$ is the neutral element of $H$. \\
For $x\in \kappa_S$, we set $d_{M,S,x}:= c_Mc_{S,x}$ where the constants $c_M$ and $c_{S,x}$ are defined in (\cite{DHSo} (1.33)). \me 

For $f\in C_c^\infty(G)$, we define the weighted orbital integral of $f$ on $G^{\si-reg}$ as follows. Let $g\in G^{\si-reg}$. We keep the above notations and  we write $g=x\ga$ with $x\in\kappa_S$ and $\ga\in S_\si$. We set
 $$\cW\cM(f)(g):=\frac{1}{d_{M,S,x}}|D_{G}(g^{-1}\si(g))|^{1/2}\int_{H\times H} f(y_1 g y_2) w_M(y_1,y_2) dy_1 dy_2.$$

 \begin{theo}\label{coefWOI}For $v,w\in V$, we have
 $$\cW\cM(c_{v,w})(g)= m_{\xi_w,\xi_v}(g),\quad g\in G^{\si-reg}.$$
 
\end{theo}
\dem Let $f_1$ be a matrix coefficient of $\tau$ and $f_2\in C_c^\infty(G)$. We set    $f:=f_1\otimes f_2$.  Let $R$ be the regular representation of $G\times G$ on $L^2(G)$ given by $[R(x_1,x_2) \Psi](g)= \Psi(x_1^{-1} gx_2)$. Then $R(f)$ is an integral operator with  smooth kernel $K_f$ given by $K_f(x,y)= \int_G f_1(xu) f_2(uy) du$. As in (\cite{DHSo} \textsection 2.2), we introduce the truncated kernel 
$$K^T(f):=\int_{H\times H} K_f(x,y) u(x,T) u(y,T) dx dy$$
where $u(x,T)$ is the truncated function of J. Arhur on $H$ (cf. \cite{DHSo} (2.7)). It is the characteristic function of a compact subset of $H$, depending on a parameter $T$ in a finite dimensional vector space, which converges to the function equal to $1$ when $\Vert T\Vert$ approaches $+\infty$. We will give the spectral asymptotic expansion of $K^T(f)$.\\
For $x\in G$, we define $$h(g):=\int_G f_1(xu) f_2(ugx) du,$$ so that $$K_f(x,y)=\big[\rho (yx^{-1})h\big](e),$$
where $\rho$ is the right regular representation of $G$.\\
If $\pi$ is a unitary irreducible admissible representation of $G$, one has 
$$\pi\big(\rho (yx^{-1})h\big)=\int_{G\times G} f_1(xu) f_2(ugy)\pi(g) du dg$$
$$=\int_{G\times G}f_1(xu) f_2(u_2) \pi(u^{-1}u_2y^{-1} )du du_2=\int_{G\times G} f_1(u_1^{-1}) f_2(u_2) \pi( u_1 x u_2 y^{-1})du_1 du_2$$
$$= \pi(\check{f_1})\pi(x)\pi(f_2) \pi(y^{-1}).$$
Since $\tau$ is supercuspidal and $f_1$ is a matrix coefficient of $\tau$, we deduce that  $\pi\big(\rho (yx^{-1})h\big)$ is equal to $0$ if $\pi$ is not equivalent to $\tau$. Therefore, applying the Plancherel formula (\cite{W2} Théorème VIII.1.1.) to $\big[\rho (yx^{-1})h]^{\check{}}$, we obtain
$$K_f(x,y)=d(\tau) \textrm{tr} \big( \tau(\check{f_1})\tau(x) \tau(f_2) \tau(y^{-1})\big).$$

We identify $ \check{V}\otimes V$ with a subspace of Hilbert-Schmidt operators on $V$. Taking an orthonormal basis $\cB_{HS}(V)$ of $ \check{V}\otimes V$ for  the scalar product $(S,S'):=\textrm{tr} (SS'^*)$, one obtains
 $$K_f(x,y)=d(\tau) \textrm{tr}  \Big( \tau(\check{f_1})\tau(x) \tau(f_2) \tau(y)^*\Big)=d(\tau)( \tau(\check{f_1})\tau(x) \tau(f_2) ,\tau(y))$$
$$=d(\tau)\sum_{S\in\cB_{HS}(V)}(\tau(\check{f_1})\tau(x) \tau(f_2) ,S^*)\overline{(\tau(y),S^*)}$$
$$=d(\tau)\sum_{S\in\cB_{HS}(V)}\textrm{tr}   \big(\tau(x) \tau(f_2)S \tau(\check{f_1})\big)\textrm{tr}   \overline{\big(\tau(y)S)},$$
where the sums over $S$ are finite since $\tau(f_2)$ and $\tau(\check{f}_1)$ are of finite rank.
Therefore, the truncated kernel is equal to
$$K^T(f)=d(\tau)\sum_{S\in\cB_{HS}(V)}P^T_{\tau}(\check{\tau}\otimes\tau(f)S) \overline{P^T_{\tau} (S)}$$
where 
$$P^T_{\tau} (S)= \int_H \textrm{tr} \big(\tau(h)S\big) u(h,T) dh, \quad S\in\check{V}\otimes V.$$
For $\check{v}\otimes v\in \check{V}\otimes V$, one has  $\textrm{tr} \big(\tau(h)(\check{v}\otimes v)\big)=c_{\check{v},v}(h)$. Since $c_{\check{v},v}$  is compactly supported, the truncated local period $P^T_{\tau} (S)$ converges when  $\Vert T\Vert$ approaches infinity to 
$$P_{\tau} (S)= \int_H \textrm{tr} \big(\tau(h)S\big)dh.$$

Therefore, we obtain
\beq\label{limspect}\lim_{\Vert T\Vert\to+\infty} K^T(f)=d(\tau)m_{P_\tau,P_\tau}(f),\eeq
where $m_{P_\tau,P_\tau}$ is the spherical character of the representation $\check{\tau}\otimes\tau$ associated to the $H\times H$-invariant linear form $P_\tau$ on $ \check{V}\otimes V$.\me

By (\cite{DHSo} Theorem  2.15), the truncated kernel $K^T(f)$ is asymptotic to a distribution $J^T(f)$ as $\Vert T\Vert$ approaches $+\infty$ and the constant term $\tilde{J}(f)$   of $J^T(f)$ is explicitly given in (\cite{DHSo} Corollary 2.11). Therefore, we deduce that  \beq\label{result1}d(\tau)m_{P_\tau,P_\tau}(f)=\tilde{J}(f).\eeq\\
We now express $m_{P_\tau,P_\tau}$ in terms of $H$-invariant linear forms on $V$.   
%%%%%%%%%%%%%%%%%%%%%%%
Let $V_H$ be the orthogonal of $V^{*H}$ in $V$. Since $\xi_u(v)=\overline{\xi_v(u)}$ for $u,v\in V$, the space $\overline{V}_H$ is the kernel of $v\mapsto \xi_v$. Let $W$ be a complementary subspace of $V_H$ in $V$. Then, the map $v\mapsto \xi_v$ is an isomorphism from $\overline{W}$ to $V^{*H}$ and $(u,v)\mapsto \xi_v(u)$ is a nondegenerate hermitian form on $W$. Let $(e_1, \ldots ,e_n)$ be an orthogonal basis of $W$ for this hermitian form. We set $\xi_i:=\xi_{e_i}$ for $i=1,\ldots, n$. Thus we have $\xi_i(e_i)\neq 0$. \\
 We identify  $\overline{V}$ and $\check{V}$ by the isomorphism $\iota$. 
We claim  that 
\beq P_\tau= \sum_{i=1}^n \frac{1}{\xi_i(e_i)}\overline{\xi_i}\otimes \xi_i\eeq
%\frac{((\xi_i,\xi_i))}{\xi_i(e_i)}
 Indeed, we have $P_\tau(v\otimes u)=\xi_v(u)=\overline{\xi_u(v)}$. Hence, the two sides are equal  to $0$ on $\overline{V}\otimes V_H+\overline{V}_H\otimes V+\overline{V}_H\otimes V_H$ and take the same value $ \xi_k(e_l)$ on $e_k\otimes e_l$ for $k,l\in\{1,\ldots n\}$. 
 %%%%%%%%%%%%%%%%%%%%%%%%%%%%%%%%%%%%%%%
Hence, by definition  of spherical characters, we deduce that   
$$m_{P_\tau,P_\tau}(f_1\otimes f_2)=\sum_{u\otimes v\in \; o.b.(\bar{V}\otimes V)}P_\tau\Big(\bar{\tau}(f_1)\otimes \tau(f_2)(u\otimes v)\Big) \overline{P_\tau(u\otimes v)}$$
$$=\sum_{u\otimes v\in \; o.b.(\bar{V}\otimes V)}\sum_{i,j=1}^n\frac{1}{\xi_i(e_i)\xi_j(e_j)}\overline{\xi_i}(\bar{\tau}(f_1)u)\xi_i(\tau(f_2)v)\overline{\overline{\xi_j}(u)\xi_j(v)},$$
where $o.b.(\bar{V}\otimes V)$ is an orthonormal basis of $\bar{V}\otimes V$.
By definition of $\bar{\xi}$ for $\xi\in V^{*H}$, one has
 $\bar{\xi}(\bar{\tau}(f_1)u)= \overline{\xi(\tau(\bar{f_1}))}.$ 
Therefore, we obtain 
\beq\label{Ptau}m_{P_\tau,P_\tau}(f_1\otimes f_2)=\sum_{i,j=1}^n\frac{1}{\xi_i(e_i)\xi_j(e_j)}\overline{m_{\xi_i,\xi_j}(\bar{f_1})}m_{\xi_i,\xi_j}(f_2).\eeq
Let $v$ and $w$ in $V$. Let $f_1:=c_{v,w}$ so that $\bar{f_1}=\check{c}_{v,w}$. . If $v\in V_H$ or $w\in V_H$, it follows from   Lemma \ref{calculcoef} that $m_{\xi_i,\xi_j}(\bar{f_1})=0$ for $i,j\in\{1,\ldots, n\}$, hence 
 $m_{P_\tau,P_\tau}(f_1\otimes f_2)=0$. Thus, we deduce from  (\ref{result1}) that
\beq\label{VH}\tilde{J}(c_{v,w}\otimes f_2)=0,\quad  v\in V_H \textrm{ or } w\in V_H.\eeq
Let $k,l\in\{1,\ldots ,n\}$. We set  $f_1:=c_{e_k,e_l}$, hence $\bar{f_1}= \check{c}_{e_l,e_k}$. By Lemma \ref{calculcoef},   one has 
$m_{\xi_i,\xi_j}( \bar{f_1})=d(\tau)^{-1} \xi_i(e_l) \xi_j(e_k).$ Therefore,  by  (\ref{result1}) and (\ref{Ptau})  we obtain 
\beq\label{ekl}\tilde{J}(c_{e_k,e_l}\otimes f_2)=m_{\xi_l,\xi_k}(f_2).\eeq
By sesquilinearity, ones deduces from (\ref{VH}) and (\ref{ekl}) that one has
\beq\label{V} \tilde{J}(c_{v,w}\otimes f_2)= m_{\xi_w, \xi_v}(f_2) \quad v,w\in V.\eeq
Let $g\in G^{\si-reg}$. Let $(J_n)_n$  be a sequence of compact open sugroups whose intersection is equal to the neutral element of $G$. The characteristic function $\phi_n$ of $J_ng J_n$ approaches the Dirac measure at $g$ as $n$ approaches $+\infty$.  Thus, if  $v,w\in V$ then   $m_{\xi_w,\xi_v}(\phi_n)$ converges to $m_{\xi_w,\xi_v}(g)$. By (\cite{DHSo} Corollary 2.11)
 the constant term $\tilde{J}(c_{v,w}\otimes \phi_n)$ converges to $ \cW\cM(c_{v,w})(g)$. We deduce the Theorem from (\ref{V}).\qed


\begin{thebibliography}{99}
\bibitem[Ar1]{ArWOI}  J. Arthur, {\em The characters of supercuspidal representations as weighted orbital integrals}, Proc. Indian Acad. Sci. Math. Sci., 97 (1987), 3-19.
\bibitem[Ar2]{ArLT}  J. Arthur, {\em A Local Trace formula}, Publ. Math. Inst. Hautes \'Etudes Sci. , 73 (1991), 5 - 96.
\bibitem[DHS]{DHSo} P. Delorme, P. Harinck and S. Souaifi, {\em Geometric side of a local relative trace formula}, arXiv:1506.09112 (47 p.), 
\bibitem[Ha]{Ha} J. Hakim, Admissible distributions on p-adic symmetric spaces, J. Reine Angew. Math. 455 (1994), 1–19.
 \bibitem[RR]{RR}  C. Rader,  S. Rallis, 
 {\em  Spherical characters on $p$-adic symmetric spaces},
Amer. J.  Math., Vol 118, N\textsuperscript{o} 1 (5 Feb. 1996),  91-178.
\bibitem[R]{R} D. Renard, {\em Représentations des groupes réductifs $p$-adiques}, Cours spécialisés, volume 17, SMF.
 \bibitem[W2]{W2}  J.-L. Waldspurger, 
{\em  La formule de Plancherel pour les groupes
$p$-adiques (d'après Harish-Chandra)},
J. Inst. Math. Jussieu  {\bf 2}  (2003)  235 - 333.
\bibitem[Z]{Z} C. Zhang, {\em Local periods for discrete series representations}, Preprint,   arXiv:1509.06166. 
\end{thebibliography}
\end{document}